\documentclass[hidelinks,onefignum,onetabnum]{siamart251216}

\usepackage{amsfonts}
\usepackage{amssymb}
\usepackage{graphicx}
\usepackage{algpseudocode}
\usepackage{multirow}
\usepackage{booktabs}
\usepackage{makecell}
\usepackage{xurl}
\DeclareGraphicsExtensions{.pdf}
\graphicspath{{pic-pdf/}}

\newcommand{\CFCA}{CF-UA}

\newcommand{\CFUA}{CF-UAMG}
\newcommand{\Pone}{GA}
\newcommand{\PW}{PW}
\newcommand{\CWM}{CWM}
\newcommand{\MIS}{MIS}

\crefname{ALG@line}{line}{lines}
\Crefname{ALG@line}{Line}{Lines}

\headers{Coarse/Fine Splitting-Based Parallel UA-AMG}
{Z. Fan, S. Hu, S. Huang, and X. Xu}

\title{\CFUA{}: A Coarse/Fine Splitting-Based Parallel Unsmoothed Aggregation Algebraic Multigrid Method\thanks{%
\textbf{Funding:} Supported by Science Challenge Project, No.~TZ2025007.}}

\author{Zhibo Fan\thanks{Laboratory of Computational Physics, Institute of Applied Physics and Computational Mathematics, Beijing 100094, 
China, and Graduate School of China Academy of Engineering Physics, Beijing, China (\email{fanzhibo23@gscaep.ac.cn}).}
\and Shaoliang Hu\thanks{Laboratory of Computational Physics, Institute of Applied Physics and Computational Mathematics, Beijing 100094, China, and CAEP Software Center for High Performance Numerical Simulation, Beijing, China (\protect\href{mailto:hu_shaoliaong@iapcm.ac.cn}{hu\_shaoliaong@iapcm.ac.cn}).}
\and Silu Huang\thanks{Laboratory of Computational Physics, Institute of Applied Physics and Computational Mathematics, Beijing 100094, 
China (\protect\href{mailto:huang_silu@iapcm.ac.cn}{huang\_silu@iapcm.ac.cn}).}
\and Xiaowen Xu\thanks{Corresponding author. Laboratory of Computational Physics, Institute of Applied Physics and Computational Mathematics, Beijing 100094, 
China (\email{xwxu@iapcm.ac.cn}).}}

\ifpdf
\hypersetup{
  pdftitle={CF-UAMG: A Coarse/Fine Splitting-Based Parallel Unsmoothed Aggregation Algebraic Multigrid Method},
  pdfauthor={Zhibo Fan, Shaoliang Hu, Silu Huang, and Xiaowen Xu}
}
\fi

\begin{document}
\maketitle

\begin{abstract}
In this paper, we propose a coarse/fine splitting-based parallel unsmoothed aggregation algebraic multigrid (AMG) method, referred to as \CFUA{}.
The proposed method utilizes the strong-connection structure identified by classical AMG (C-AMG) coarse/fine splitting to guide the selection of aggregate roots and the formation of aggregates, aiming to improve the quality of the aggregation-based coarse space and thus the convergence of unsmoothed aggregation AMG (UA-AMG).
We also prove uniform two-grid convergence of \CFUA{} with a convergence bound independent of the matrix size for symmetric positive definite M-matrices with nonnegative row sums.
Numerical experiments on three classes of three-dimensional problems, including Poisson-type equations, radiation diffusion equations, and incompressible Navier--Stokes equations, demonstrate the robustness both in convergence and parallelism of \CFUA{}.
Compared with UA-AMG methods employing four commonly used aggregation algorithms, \CFUA{} achieves average speedups of $1.29\times$, $2.04\times$, and $1.94\times$ for the three classes of problems, respectively.
\end{abstract}



\section{Introduction}\label{sec:intro}
In scientific and engineering computing, many numerical simulations based on partial differential equations (PDEs) typically require the solution of linear systems of the form
\begin{equation}
    \label{Ax=b}
    A x = b,
\end{equation}
where $A\in\mathbb{R}^{n\times n}$ is a sparse matrix, and $b,x\in\mathbb{R}^{n}$ are the right-hand side and solution vectors, respectively.
Solving such systems has become a major computational bottleneck in many real-world applications as the scale of numerical simulations increases.
Preconditioned Krylov subspace methods are widely used to solve such systems, and their performance depends critically on the construction of preconditioners~\cite{2003-saad-book}.
Algebraic multigrid (AMG)~\cite{1999-AMG} is one of the most widely used preconditioners for large-scale problems owing to its optimal or near-optimal computational complexity~\cite{XXW-2019parallel}.

AMG methods can be classified into two main categories, classical AMG (C-AMG) and aggregation-based AMG, according to their coarsening strategies.
C-AMG partitions the set of degrees of freedom $\Omega$ into a coarse subset $\Omega_C$ and a fine subset $\Omega_F$ through coarse/fine (C/F) splitting, with $\Omega_C$ defining the coarse grid~\cite{1999-AMG}.
Aggregation-based AMG, in contrast, groups the degrees of freedom into disjoint sets, called aggregates, and combines the variables within each aggregate into a single coarse variable, with the resulting coarse variables defining the coarse grid~\cite{1996-SAMG,2005-Adaptive-SA}.
Aggregation-based AMG employs a simpler coarsening procedure and typically constructs an initial piecewise constant interpolation operator from the aggregates, leading to low operator complexity.
Moreover, aggregation-based coarsening and interpolation can preserve certain discrete conservation properties on coarse levels~\cite{1995-N-UA-AMG,1996-N-UA-AMG}, which contributes to improved robustness for nonsymmetric and coupled problems, such as convection--diffusion and Navier--Stokes problems~\cite{2006-UA-PW,2008-Muresan-Notay-UA,2008-Kcycle,2012-UA-CD,2016-Notay-Stokes,2017-Notay-Stokes,2022-Notay-NS,2025-Notay-NS}.
Aggregation-based AMG methods can be further classified into two categories, smoothed aggregation AMG (SA-AMG) and unsmoothed aggregation AMG (UA-AMG), depending on whether the interpolation operator is smoothed.
SA-AMG enhances the approximation properties of the interpolation operator through smoothing, improving convergence at the expense of increased operator complexity.
UA-AMG, by contrast, is widely used in large-scale applications owing to its low operator complexity~\cite{2008-Muresan-Notay-UA,2012-UA-CD}.

The efficiency of UA-AMG depends critically on aggregate quality.
The classical greedy aggregation (\Pone{}) algorithm~\cite{1996-SAMG,1995-N-UA-AMG,1996-N-UA-AMG} constructs aggregates using a local greedy strategy, sequentially selecting unaggregated variables as roots and grouping each root with its unaggregated neighbors in the matrix adjacency graph.
This algorithm typically produces relatively large aggregates and correspondingly small coarse systems.
However, the resulting aggregation is strongly influenced by the order in which roots are selected, and root selection does not incorporate global information.
To address this limitation, independent-set-based aggregation (\MIS{}) algorithms~\cite{2000-par-sa,2013-GPU-ParUA} have been developed.
Specifically, the \MIS{} aggregation algorithm selects a maximal set of vertices whose pairwise graph distances are greater than two and uses these vertices as aggregate roots, which makes the algorithm suitable for parallel implementation.
However, it does not fully exploit the strength-of-connection structure of the matrix, which can limit the convergence of the resulting AMG method.
By comparison, the pairwise (\PW{}) aggregation algorithm~\cite{2006-UA-PW,2012-UA-CD} forms relatively small aggregates by pairing variables, yielding correspondingly larger coarse systems while effectively preserving the strong local couplings.
Further approaches have been proposed to improve aggregate quality by incorporating additional information into these basic aggregation algorithms~\cite{2005-Adaptive-SA,2013-CWM-UA,2021-parallel-CWM,2018-coordinate-aggregation,2022-EdgeMatrix-AMG}.
For example, the compatible weighted matching (\CWM{}) aggregation algorithm~\cite{2013-CWM-UA,2021-parallel-CWM} relies on a representative test vector to guide aggregate formation, although obtaining such a vector by relaxation can be challenging for complex problems.
Geometric information has also been used to improve aggregate construction~\cite{2018-coordinate-aggregation,2022-EdgeMatrix-AMG}.
Despite improvements in aggregate quality, these approaches rely on additional information that is either extracted from the underlying problem or obtained through additional computation.

We note that although both C-AMG and UA-AMG utilize the local information in the matrix adjacency graph, the ways of utilization are essentially different.
More specifically, C-AMG constructs the coarse grid through C/F splitting, using the overall structure of the strength-of-connection graph to guide coarse-point selection.
In contrast, UA-AMG selects aggregate roots and forms aggregates primarily using local connectivity, and the quality of its coarse space depends strongly on the representativeness of the selected roots.
Motivated by this observation, we develop a hybrid coarsening strategy that combines C/F splitting with aggregation.
The C-points selected by C/F splitting serve as initial aggregate roots, and strong connections between C-points and F-points guide the assignment of F-points to aggregates.
This strategy incorporates the global structure information captured by C/F splitting directly into UA-AMG aggregation to improve the quality of the coarse space without requiring additional information.
We denote the proposed parallel UA-AMG method by \CFUA{}.

The remainder of this paper is organized as follows.
The basics of C-AMG and UA-AMG are reviewed in \cref{sec:AMG-base}.
The proposed aggregation algorithm and its parallel implementation are described in \cref{sec:aggregation}.
A uniform two-grid convergence bound for \CFUA{} is established in \cref{sec:convergence}.
Numerical results for three representative classes of problems are presented in \cref{sec:numerics}.
Conclusions are drawn in \cref{sec:conclusion}.

\section{Algebraic multigrid basics}
\label{sec:AMG-base}
\subsection{Classical AMG}
Consider the linear system~\eqref{Ax=b}.
The set of degrees of freedom is denoted by $\Omega = \{1,2,\ldots,n\}$.
C-AMG consists of a setup phase and a solve phase, as described for the two-grid case in \cref{alg:two-level-amg}~\cite{1999-AMG}.
The setup phase determines the coarse grid and constructs the interpolation operator $P$, the restriction operator $R$, and the coarse-grid matrix $A_c$.
The solve phase consists of two processes: smoothing, typically performed using Jacobi or Gauss--Seidel iterations, and coarse-grid correction (CGC), which reduce high- and low-frequency error components, respectively.
\begin{algorithm}[htbp]
    \caption{Two-grid C-AMG algorithm.}
\label{alg:two-level-amg}
\begin{algorithmic}[1]
        \Statex $\blacktriangleright$~\textbf{Setup phase:}
        \State \label{MG_step:CF}Perform C/F splitting $\Omega = \Omega_C \cup \Omega_F$ and select $\Omega_C$ as the coarse grid.
        \State Construct interpolation $P:\mathbb{R}^{n_c}\rightarrow\mathbb{R}^{n}$ and restriction $R:\mathbb{R}^{n}\rightarrow\mathbb{R}^{n_c}$, where $n_c=|\Omega_C|$. Compute $A_c=RAP$.
        \Statex $\blacktriangleright$~\textbf{Solve phase:}
        \State Apply $\nu_1$ presmoothing steps: $x \leftarrow \text{Pre-smoothing}(A,b,x,\nu_1)$.
        \State Compute the residual and restrict it to the coarse grid: $r_c \leftarrow R ( b - Ax)$.
        \State \label{MG_step:CGC}Solve the coarse-grid correction equation: $A_c \delta_c = r_c$.
        \State Interpolate the error to the fine grid and correct the approximate solution: $x \leftarrow x + P \delta_c$.
        \State Apply $\nu_2$ postsmoothing steps: $x \leftarrow \text{Post-smoothing}(A,b,x,\nu_2)$.
    \end{algorithmic}
\end{algorithm}

A multilevel AMG method is obtained by recursively applying the two-grid procedure to the coarse-grid system in step~\ref{MG_step:CGC} when the system remains large and solving it exactly is expensive.
Standard multilevel cycling strategies include V- and W-cycles.

\subsection{Unsmoothed aggregation AMG}\label{sub-section: UA-AMG}
\begin{figure}[htbp]
    \centering
    \includegraphics[width = 1.0\textwidth]{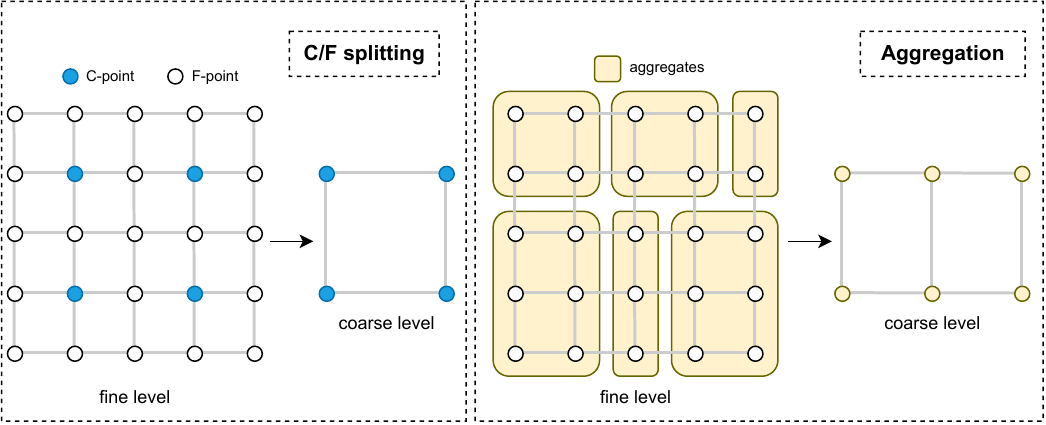}
    \caption{Comparison of coarsening in C-AMG (left) and UA-AMG (right).}
    \label{fig:cf-vs-agg}
\end{figure}

UA-AMG and C-AMG differ primarily in their coarsening procedures.
C-AMG defines the coarse grid $\Omega_C$ through C/F splitting of $\Omega$ into coarse and fine subsets (step~\ref{MG_step:CF} of \cref{alg:two-level-amg}).
UA-AMG, by contrast, constructs the coarse grid by partitioning $\Omega$ into aggregates
\begin{equation}
    \mathcal{G} = \left\{  G_i  \right\}_{i=1}^{n_c},
\end{equation}
where the nonempty sets $G_i\subset\Omega$ satisfy $G_i\cap G_j=\emptyset$ for $i\neq j$ and $\bigcup_{i=1}^{n_c}G_i=\Omega$.
The variables in each aggregate are combined into a single coarse variable, yielding a coarse space of dimension $n_c$.
\Cref{fig:cf-vs-agg} compares the coarsening procedures of C-AMG and UA-AMG\@.

The piecewise constant interpolation operator associated with $\mathcal{G}$ is defined by
\begin{equation}
    (P)_{ij} =
    \begin{cases}
        1, & i \in G_j, \\
        0, & i \notin G_j.
    \end{cases}
\end{equation}
The Galerkin coarse-grid matrix $A_c=RAP$, with $R=P^T$, is assembled by summing the matrix entries over pairs of aggregates:
\begin{equation}
    (A_c)_{ij}
    =
    \sum_{k\in G_i}
    \sum_{l\in G_j}
    a_{kl},
    \qquad
    1\leq i,j\leq n_c.
\end{equation}

\section{C/F splitting-based parallel aggregation}\label{sec:aggregation}
Existing aggregation algorithms have issues with convergence and parallel scalability in real-world applications.
Regarding convergence, existing aggregation algorithms often give insufficient consideration to the global approximation properties of the resulting coarse space.
They typically construct aggregates using local connection strengths in the matrix adjacency graph.
In addition, some algorithms can produce different aggregates depending on the order in which variables are processed, leading to variations in iteration counts as the number of processes increases.
Convergence can be improved by incorporating smooth-error information. However, obtaining sufficiently accurate information remains challenging for complex problems.
Regarding parallel scalability, some algorithms reduce communication overhead by constructing aggregates independently within each process.
This approach simplifies parallel implementation at the cost of making the resulting aggregate structure dependent on the domain partitioning.
An alternative approach is to construct aggregates across process boundaries using interprocess coupling information, typically at the cost of more complex communication and load-balancing adjustments.
The associated communication overhead can increase with the number of processes.

\begin{figure}[htbp]
    \centering
    \includegraphics[width = 0.8\textwidth]{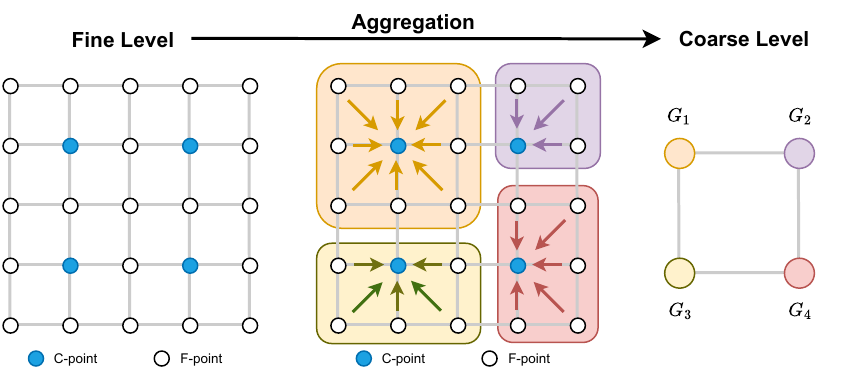}
    \caption{Overview of the C/F splitting-based aggregation algorithm (\CFCA{}).}
    \label{CF-Agg-framework}
\end{figure}

We propose a C/F splitting-based parallel aggregation algorithm to address these issues.
C-AMG uses the overall structure of the strength-of-connection graph in C/F splitting to select representative coarse points with a balanced distribution~\cite{xu2017-AI-Coarse}.
In contrast, UA-AMG selects aggregate roots and forms aggregates primarily using local connectivity, and the quality of its coarse space depends strongly on the representativeness of these roots.
We therefore use the C-points selected by C/F splitting as initial aggregate roots and guide the assignment of F-points to aggregates using strong connections.
This approach incorporates global structure information into aggregation to improve the quality of the aggregation-based coarse space.
The resulting parallel UA-AMG method is denoted by \CFUA{}, and its aggregation algorithm by \CFCA{}.
The construction is illustrated in \cref{CF-Agg-framework} and described below.

\subsection{Aggregation procedure}
The set of degrees of freedom $\Omega$ is partitioned into $p$ disjoint subsets $\Omega_q$, one for each process, with $\Omega=\bigcup_{q=1}^{p}\Omega_q$ and $\Omega_q\cap\Omega_r=\emptyset$ for $q\neq r$.
The \CFCA{} algorithm is given in \cref{C/F-Agg}.

\begin{algorithm}[htbp]
    \caption{C/F splitting-based parallel aggregation algorithm (\CFCA{}).}
\label{C/F-Agg}
\begin{algorithmic}[1]
        \Require Matrix $A=(a_{ij})\in\mathbb{R}^{n\times n}$ and the set of degrees of freedom $\Omega$.
        \Ensure Collection of aggregates $\mathcal{G}=\{G_k\}_{k=1}^{n_c}$.
        \Statex $\blacktriangleright$~\textbf{Stage 1: Partitioning and root initialization}
        \State For each process $q$, partition $\Omega_q$ into the interior point set $\Omega^I_q$ and the boundary point set $\Omega^B_q$.
        \State Apply parallel C-AMG coarsening to obtain the C/F splitting $\Omega = \Omega_C \cup \Omega_F$, and initialize the collection of aggregates using the points in $\Omega_C$ as roots.
        \Statex $\blacktriangleright$~\textbf{Stage 2: Boundary reclassification}
        \State For each point $i\in\Omega^B_q$, if $i$ is an endpoint of a strong cross-process connection involving a C-point, retain $i$ in $\Omega^B_q$. Otherwise, move $i$ from $\Omega^B_q$ to $\Omega^I_q$. The reclassification involves three situations, as detailed in \cref{sec:CFUA-Stage2}.
        \Statex $\blacktriangleright$~\textbf{Stage 3: Aggregation of interior and boundary points}
        \State \textbf{Interior points:} For each point $i\in\Omega^I_q$, identify the C-point $j^*$ on which $i$ depends most strongly. Three aggregation situations are distinguished according to the strong-connection relation between $i$ and $j^*$, as detailed in \cref{sec:CFUA-Stage3}.
        \State \label{CFUA-algorithm:boundary}\textbf{Boundary points:} Aggregate the points in $\Omega^B_q$ according to the same rules used for the interior points in $\Omega^I_q$.
    \end{algorithmic}
\end{algorithm}

The three stages of \cref{C/F-Agg} are described in detail below.

\subsubsection{Stage 1: Partitioning and root initialization}\label{sec:CFUA-Stage1}
The local set $\Omega_q$ is partitioned into interior points $\Omega^I_q$ and boundary points $\Omega^B_q$ according to the matrix adjacency structure:
\begin{equation*}
    \Omega^I_q = \{
        i\in\Omega_q \mid a_{ij} = 0 \text{ for all } j \notin \Omega_q
    \}, \qquad
    \Omega^B_q = \Omega_q \setminus \Omega^I_q.
\end{equation*}
C/F splitting $\Omega=\Omega_C\cup\Omega_F$ is obtained using a parallel coarsening strategy in C-AMG.
The C-points serve as the initial aggregate roots, defining the initial collection of aggregates $\mathcal{G}$ and the aggregation map $m(\cdot)$:
\begin{equation*}
     \mathcal{G} = \{\{i\}:
     i\in \Omega_C \}, \qquad m(i) = \text{aggregate\_id}(i),
\end{equation*}
where $m(i)$ is the index of the aggregate containing variable $i$.

PMIS~\cite{2006-PMIS-HMIS} is an inherently parallel coarsening algorithm that is widely used in real-world applications~\cite{XXW-2019parallel} and typically yields low operator complexity.
We therefore employ PMIS as the parallel C/F splitting algorithm in \CFCA{}, consistent with the low-complexity characteristics of UA-AMG\@.

\subsubsection{Stage 2: Boundary reclassification}\label{sec:CFUA-Stage2}
Boundary aggregation in Stage~3 involves searching for C-points on other processes (step~\ref{CFUA-algorithm:boundary} of \cref{C/F-Agg}).
If the C-points on which a boundary F-point strongly depends are all local, these searches cannot find a suitable C-point and incur unnecessary communication.
We therefore reclassify such boundary F-points as interior points before Stage~3, allowing them to be aggregated locally without these searches or the associated communication.
Specifically, both endpoints of a strong cross-process connection involving a C-point are retained as boundary points, and all other boundary points are reclassified as interior points.
The retained boundary points include both F-points requiring nonlocal root information and C-points supplying that information.

\begin{figure}[htbp]
    \centering
    \includegraphics[width = 1.0\textwidth]{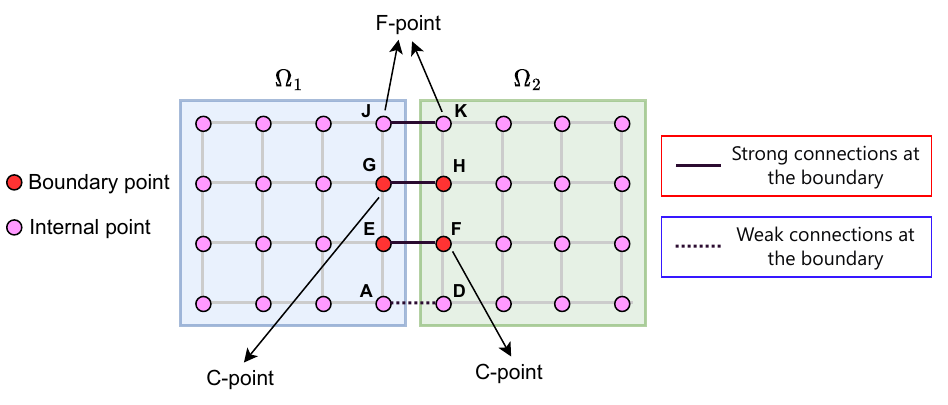}
    \caption{Illustration of boundary point reclassification for two processes.}
    \label{fig:cfb-partition}
\end{figure}

\Cref{fig:cfb-partition} illustrates the partition into interior and boundary points and the boundary reclassification in Stage~2 using a two-process example, where $\Omega_1$ and $\Omega_2$ denote the subsets of $\Omega$ assigned to the two processes.
The eight boundary points are \textbf{A}, \textbf{D}, \textbf{E}, \textbf{F}, \textbf{G}, \textbf{H}, \textbf{J}, and \textbf{K}.
Points \textbf{G} and \textbf{F} are C-points, whereas \textbf{J} and \textbf{K} are F-points.
The connections \textbf{E}--\textbf{F}, \textbf{G}--\textbf{H}, and \textbf{J}--\textbf{K} are strong (solid lines), whereas the connection \textbf{A}--\textbf{D} is weak (dashed line).
For boundary reclassification, we distinguish the following three situations based on the connection strengths and the C/F classifications of the endpoints:
\begin{itemize}
    \item \textbf{Case 1:} The weak cross-process connection \textbf{A}--\textbf{D} does not contribute to aggregate assignment, so both points are reclassified as interior points.
    \item \textbf{Case 2:} The strong cross-process connections \textbf{E}--\textbf{F} and \textbf{G}--\textbf{H} each involve a C-point. Both endpoints of each connection are retained as boundary points.
    \item \textbf{Case 3:} The cross-process connection \textbf{J}--\textbf{K} involves only F-points. Both are reclassified as interior points, despite their strong connection.
\end{itemize}

\subsubsection{Stage 3: Aggregation of interior and boundary points}\label{sec:CFUA-Stage3}

\begin{figure}[htbp]
    \centering
    \includegraphics[width = 1.0\textwidth]{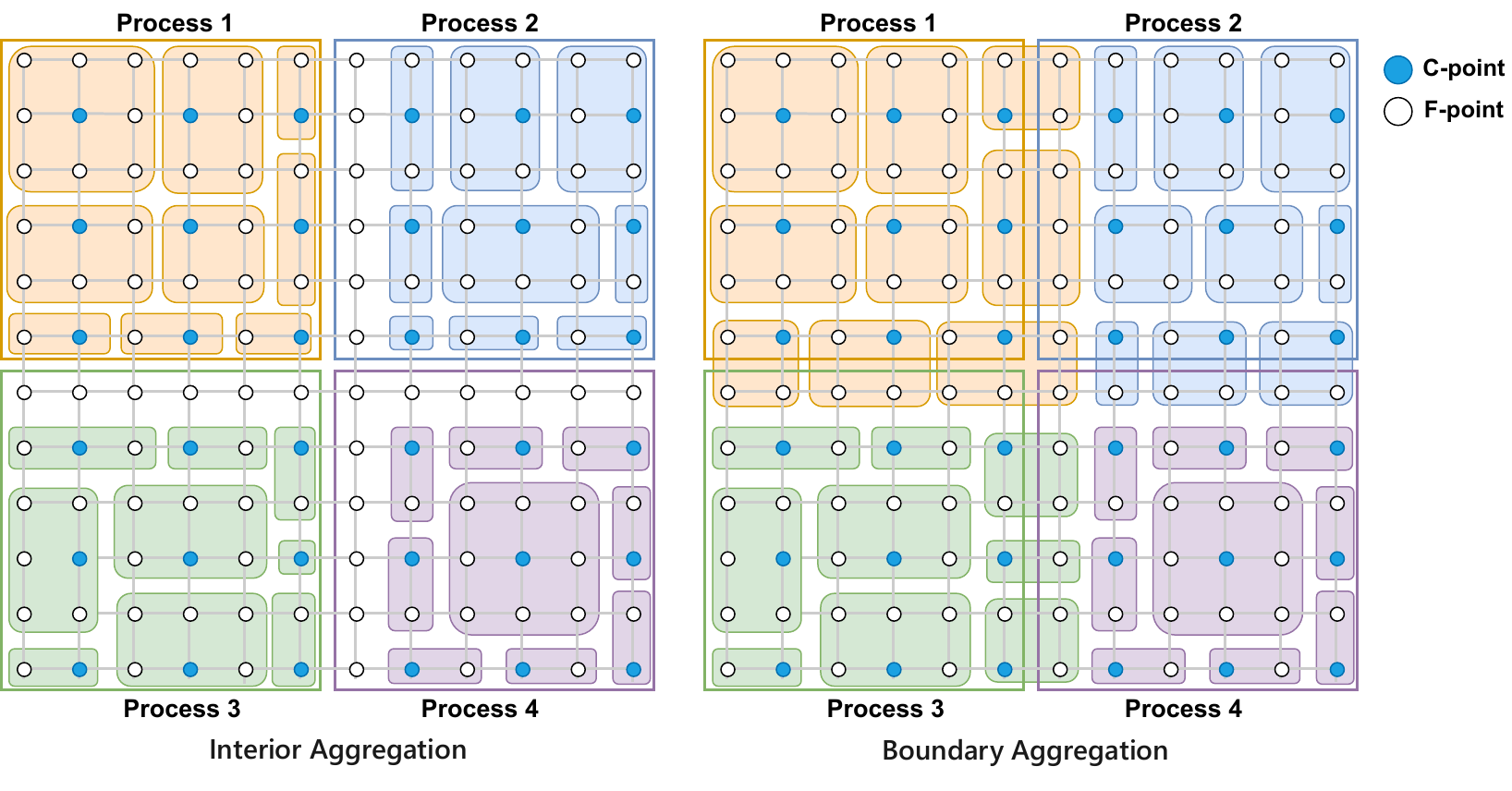}
    \caption{Illustration of the aggregation of interior and boundary points.}
    \label{CF-Agg-framework-Par}
\end{figure}

Each process first aggregates its interior points independently, as illustrated on the left of \cref{CF-Agg-framework-Par}.
For each unaggregated interior point $i\in\Omega^I_q$, $j^*=\arg\max_{k\in\Omega_C}|a_{ik}|$ is the C-point most strongly connected to $i$.
The aggregate assignment for $i$ is then determined according to the following three situations:
\begin{itemize}
    \item \textbf{Case 1:} If $i$ strongly depends on $j^*$, add $i$ to the aggregate containing $j^*$.
    \item \textbf{Case 2:} If $i$ does not strongly depend on any C-point, search among the local non-C-point roots on which $i$ strongly depends. Add $i$ to the aggregate of the root on which it depends most strongly.
    \item \textbf{Case 3:} If neither Case~1 nor Case~2 applies, create a new aggregate with $i$ as its root.
\end{itemize}

The processes synchronize the root aggregate information after aggregating the interior points.
The remaining boundary points are then assigned using the same rules, as illustrated on the right of \cref{CF-Agg-framework-Par}.

\subsection{Algorithmic properties}
\label{sec:cfua-features}
The \CFUA{} method improves the convergence robustness and parallel scalability of UA-AMG by incorporating the C/F splitting of C-AMG into aggregation-based coarsening, while maintaining the low operator complexity of UA-AMG\@.

Regarding convergence robustness, the \CFCA{} algorithm assigns every point to an aggregate and ensures that each nonroot point strongly depends on the root.
This property is central to the uniform two-grid convergence analysis in \cref{sec:convergence}.
Moreover, C/F splitting provides representative initial roots and reduces the sensitivity of the resulting aggregates to the order in which roots are selected.
Strong C--F connections guide the assignment of F-points to aggregates, and the global structure information supplied by C/F splitting improves the quality of the aggregation-based coarse space without requiring additional information.
The aggregation rules in Stage~3 also handle points without a suitable C-point root, ensuring that no point is left unaggregated.

Regarding parallel scalability, the \CFCA{} algorithm has low communication overhead and enables overlap between communication and computation.
After C/F splitting, aggregate membership at process boundaries is largely determined, reducing the need for dynamic reassignment.
Boundary reclassification further reduces unnecessary communication by retaining only the endpoints of strong cross-process connections involving C-points as boundary points.
Apart from the communication required for C/F splitting, the algorithm requires only one exchange of root aggregate information in Stage~3.
The remaining aggregation operations are performed locally and in parallel across processes.
Since interior points can be aggregated using only local information, this communication can overlap with their aggregation to improve parallel efficiency.

\section{Two-grid convergence of \CFUA{}}
\label{sec:convergence}
We establish uniform two-grid convergence of \CFUA{} for symmetric positive definite M-matrices with nonnegative row sums.
Let $\Pi = P A_c^{-1} P^T A$ denote the coarse-grid projection for \CFUA{}, and let $T = I - \Pi$ denote the corresponding coarse-grid correction operator.
Let $G$ denote the error propagation matrix of the smoother.
The two-grid error propagation operator with postsmoothing is $E_{\mathrm{TG}}=GT$.

We use the theoretical results in~\cite{maclachlan2014theoretical} and first recall the following theorem.
\begin{theorem}\label{thm:MO}
Let $A$ be symmetric positive definite, and let $g$ be a positive-definite function, with $g(0)=0$ and $g(e)>0$ for every nonzero $e\in\mathbb{R}^n$.
Suppose that there exist constants $\bar{\alpha}_g, \bar{\beta}_g > 0$ such that, for every $e$,
\begin{align}
    \|G e\|_A^2 &\leq \|e\|_A^2 - \bar{\alpha}_g \, g(e)
    &&\text{(smoothing property)}, \label{eq:MO-smooth}\\
    \|T e\|_A^2 &\leq \bar{\beta}_g \, g(T e)
    &&\text{(approximation property)}. \label{eq:MO-approx}
\end{align}
Then
\begin{equation*}
    \|E_{\mathrm{TG}}\|_A \leq \sqrt{1 - \frac{\bar{\alpha}_g}{\bar{\beta}_g}}.
\end{equation*}
\end{theorem}

The proof is given in~\cite{maclachlan2014theoretical}.
The convergence estimate in \cref{thm:MO} relies on the smoothing property~\eqref{eq:MO-smooth} and the approximation property~\eqref{eq:MO-approx}.
We choose $g(e):=\|Ae\|_{D^{-1}}^2$, where $D=\operatorname{diag}(A)$.
For the smoothing property, we invoke standard estimates for smoothers such as $\omega$-Jacobi and Gauss--Seidel and assume that $\bar{\alpha}_g>0$ can be chosen independently of the matrix size $n$~\cite{2003-saad-book,1999-AMG}.
We therefore focus on establishing the approximation property for \CFUA{}.

\subsection{Approximation property}
Let $d_i=a_{ii}$, $i=1,\ldots,n$, denote the diagonal entries, and set $w_{ij}=-a_{ij}=|a_{ij}|$ for $i\neq j$.
With the row sums denoted by $\sigma_i=d_i-\sum_{j\neq i}w_{ij}$, the energy admits the decomposition~\cite{1999-AMG}
\begin{equation}
    \|x\|^2_A = \sum_{i<j} w_{ij}(x_i - x_j)^2 + \sum_i \sigma_i x_i^2.
\end{equation}

For an aggregate $G\in\mathcal{G}$ and any vector $x$, define the local energy by
\begin{equation}
    E_G(x) = \sum_{i,j\in G,\, i<j} w_{ij}(x_i - x_j)^2.
\end{equation}
The energy decomposition then implies
\begin{equation}
    \sum_{G\in\mathcal{G}} \left[E_G(x) + \sum_{i\in G} \sigma_i x_i^2\right]
    \leq \|x\|^2_A.
\end{equation}
Equality holds when there are no connections between distinct aggregates.

We also define the sum of the diagonal entries over $G$ by
\begin{equation}
    s_G := \sum_{i\in G} d_i.
\end{equation}
The weighted mean of $x$ on $G$ is
\begin{equation}
    \bar{x}_G := \frac{1}{s_G} \sum_{i\in G} d_i x_i
\end{equation}
and the corresponding weighted variance is
\begin{equation}
    V_G(x) := \sum_{i\in G} d_i(x_i - \bar{x}_G)^2.
\end{equation}
The following lemma establishes the coarse-grid approximation property in terms of these local quantities.

\begin{lemma}\label{lem:weak-approximation}
Let $A$ be a symmetric positive definite M-matrix.
Suppose that, for each aggregate $G\in\mathcal{G}$, there exists a constant $\mu_G > 0$ such that, for every $x\in\mathbb{R}^n$,
\begin{equation}
    \mu_G V_G(x) \leq \left[E_G(x) + \sum_{i\in G} \sigma_i x_i^2\right].
\end{equation}
Then, for every $e\in\mathbb{R}^n$,
\begin{equation}
    \|T e\|_A^2 \leq \bar{\beta}_g \|A T e\|_{D^{-1}}^2,
\end{equation}
where
\begin{equation}
    \bar{\beta}_g = \frac{1}{\min\limits_{G\in\mathcal{G}}\mu_G}.
\end{equation}
\end{lemma}

\begin{proof}
For any $e$, set
\begin{equation*}
    z = T e = e - \Pi e.
\end{equation*}
Since $\Pi$ is an $A$-orthogonal projection, $z$ is $A$-orthogonal to the range of $P$, and hence $P^TAz=0$.
By the structure of $P$, this condition is equivalent to
\begin{equation}
    \sum_{i\in G}(A z)_i = 0
\end{equation}
for every aggregate $G$.
Define $\delta_i=z_i-\bar{z}_G$ for $i\in G$, and let $G(i)$ denote the aggregate containing $i$.
The Cauchy--Schwarz inequality then yields
\begin{align*}
    \|z\|_A^2
    &= \sum_{i=1}^n (\bar{z}_{G(i)} + \delta_i)(A z)_i \\
    &= \sum_{G\in\mathcal{G}} \bar{z}_G \sum_{i\in G}(A z)_i
       + \sum_{i=1}^n \delta_i(A z)_i \\
    &= \sum_{i=1}^n \delta_i(A z)_i
     = \sum_{i=1}^n \sqrt{a_{ii}}\,\delta_i\frac{1}{\sqrt{a_{ii}}}(A z)_i \\
    &\leq \|\delta\|_D\|A z\|_{D^{-1}},
\end{align*}
where
\begin{align*}
    \|\delta\|_D^2
    &= \sum_{i=1}^n a_{ii}(z_i - \bar{z}_{G(i)})^2 \\
    &= \sum_{G\in\mathcal{G}} \sum_{i\in G} d_i(z_i - \bar{z}_G)^2
     = \sum_{G\in\mathcal{G}} V_G(z).
\end{align*}
By hypothesis,
\begin{equation*}
    \sum_{G\in\mathcal{G}} V_G(z)
    \leq \sum_{G\in\mathcal{G}} \frac{1}{\mu_G}
    \left[E_G(z) + \sum_{i\in G}\sigma_i z_i^2\right]
    \leq \bar{\beta}_g \|z\|_A^2.
\end{equation*}
Thus $\|\delta\|_D^2 \leq \bar{\beta}_g\|z\|_A^2$.
Combining this bound with the Cauchy--Schwarz estimate yields $\|z\|_A^2 \leq \bar{\beta}_g\|A z\|_{D^{-1}}^2$.
Consequently,
\begin{equation*}
    \|T e\|_A^2 = \|z\|_A^2 \leq \bar{\beta}_g\|A z\|_{D^{-1}}^2.
\end{equation*}
\end{proof}

The approximation constant $\bar{\beta}_g$ in \cref{lem:weak-approximation} is determined by the local constants $\mu_G$, which measure aggregate quality.
Larger values of $\mu_G$ yield a smaller approximation constant and, for a fixed smoothing constant, a sharper two-grid convergence bound.
Each aggregate $G$ produced by \CFCA{} can be written as
\begin{equation}
    G = \{r\} \cup G^F,
\end{equation}
where $r$ is the unique root and $G^F$ contains the nonroot points.
This aggregate structure yields the following estimate for $\mu_G$.

\begin{lemma}\label{lem:simple-approximation}
Let $A$ be a symmetric positive definite M-matrix with nonnegative row sums, and let $\mathcal{G}$ be the collection of aggregates produced by \CFCA{}.
Then, for every $x\in\mathbb{R}^n$ and every $G\in\mathcal{G}$ with $|G|\geq 2$,
\begin{equation}
    \left(\min\limits_{i\in G^F}\frac{w_{ir}}{d_i}\right)V_G(x)
    \leq \left[E_G(x) + \sum_{i\in G}\sigma_i x_i^2\right].
\end{equation}
\end{lemma}
\begin{proof}
The minimization property of the weighted variance gives
\begin{align*}
    V_G(x)
    &= \min\limits_c \sum_{i\in G} d_i(x_i - c)^2 \\
    &\leq \sum_{i\in G} d_i(x_i - x_r)^2
     = \sum_{i\in G^F} d_i(x_i - x_r)^2.
\end{align*}
On the other hand,
\begin{equation*}
    E_G(x) = \sum_{i,j\in G,\, i<j}w_{ij}(x_i - x_j)^2
    \geq \sum_{i\in G^F} w_{ir}(x_i - x_r)^2.
\end{equation*}
Combining these inequalities with the nonnegativity of the row sums yields
\begin{align*}
    \left(\min\limits_{i\in G^F}\frac{w_{ir}}{d_i}\right)V_G(x)
    &\leq \sum_{i\in G^F}
    \left(\min\limits_{i\in G^F}\frac{w_{ir}}{d_i}\right)d_i(x_i - x_r)^2 \\
    &\leq E_G(x)
     \leq \left[E_G(x) + \sum_{i\in G}\sigma_i x_i^2\right].
\end{align*}
\end{proof}

Let $\theta\in(0,1]$ denote the strength-of-connection threshold, and let $s\geq 2$ denote the maximum number of nonzero entries in any row of $A$.
For each row with a nonzero off-diagonal entry, define the diagonal dominance ratio by
\begin{equation*}
    \kappa_i=\frac{a_{ii}}{\sum_{j\neq i}|a_{ij}|}\geq 1.
\end{equation*}

We make the following assumptions:
\begin{itemize}
    \item There exists a constant $\eta$, independent of $n$, such that $\kappa_i\leq\eta$ for all such rows.
    \item $s$ is a constant independent of $n$.
\end{itemize}

The bound on $\kappa_i$ is motivated by the effectiveness of relaxation in damping error components associated with strongly diagonally dominant rows~\cite{1999-AMG}.
The assumption on $s$ is natural for PDE discretizations with a fixed local stencil, for which the maximum number of nonzero entries per row remains constant as $n$ increases.
Under these assumptions, we obtain the following theorem.

\begin{theorem}[Approximation bound]
Let $A$ be a symmetric positive definite M-matrix with nonnegative row sums.
The approximation constant for the collection of aggregates $\mathcal{G}$ produced by \CFCA{} satisfies
\begin{equation}
    \bar{\beta}_g \leq \frac{\eta(s - 2 + \theta)}{\theta}.
\end{equation}
\end{theorem}
\begin{proof}
By the construction of \CFCA{}, every nonroot point in a nonsingleton aggregate satisfies the strength criterion
\begin{equation*}
    w_{ir} \geq \theta\max\limits_{j\neq i}w_{ij},
    \qquad \text{for all } i\in G^F.
\end{equation*}
Since each row contains at most $s$ nonzero entries,
\begin{align*}
    \sum_{j\neq i}w_{ij}
    &= w_{ir} + \sum_{j\neq i,\,j\neq r}w_{ij} \\
    &\leq w_{ir} + (s - 2)\max\limits_{j\neq i}w_{ij}
     \leq w_{ir} + \frac{s - 2}{\theta}w_{ir}.
\end{align*}
It follows that
\begin{equation*}
    \frac{w_{ir}}{\sum_{j\neq i}w_{ij}} \geq \frac{\theta}{s - 2 + \theta}.
\end{equation*}
Moreover,
\begin{equation*}
    d_i = \kappa_i\sum_{j\neq i}w_{ij}
    \leq \eta\sum_{j\neq i}w_{ij}.
\end{equation*}
Combining these estimates, we obtain
\begin{equation*}
    \frac{w_{ir}}{d_i} \geq \frac{\theta}{\eta(s - 2 + \theta)},
    \qquad \text{for all } i\in G^F.
\end{equation*}
For a singleton aggregate $G=\{r\}$, we have $\bar{x}_G=x_r$, and hence $V_G(x)=E_G(x)=0$.
The nonnegativity of the row sums implies $\sigma_r\geq 0$. Thus, for every $\mu_G>0$,
\begin{equation*}
    \mu_G V_G(x) \leq E_G(x) + \sum_{i\in G}\sigma_i x_i^2.
\end{equation*}
Together with \cref{lem:simple-approximation}, the estimates for singleton and nonsingleton aggregates allow the choice
\begin{equation*}
    \mu_G = \frac{\theta}{\eta(s-2+\theta)}
\end{equation*}
uniformly for all aggregates.
Applying \cref{lem:weak-approximation} gives the stated upper bound on the approximation constant.
\end{proof}

\subsection{Uniform two-grid convergence}
Combining the approximation bound with a uniform smoothing estimate yields the following uniform two-grid convergence result.

\begin{theorem}[Uniform two-grid convergence]\label{lem:convergence}
Let $A$ be a symmetric positive definite M-matrix with nonnegative row sums.
Suppose that the chosen smoother satisfies the smoothing property~\eqref{eq:MO-smooth} with a constant $\bar{\alpha}_g>0$.
Then the two-grid operator $E_{\mathrm{TG}}$ of \CFUA{} satisfies
\begin{equation}
\label{eq: TW-convergence}
    \|E_{\mathrm{TG}}\|_A
    \leq \sqrt{1 - \frac{\bar{\alpha}_g}{\bar{\beta}_g}},
    \qquad 0<\bar{\beta}_g \leq \frac{\eta(s - 2 + \theta)}{\theta},
\end{equation}
where $\bar{\beta}_g$ is an approximation constant satisfying~\eqref{eq:MO-approx}.
If $\theta$ is fixed, $\eta$ and $s$ are uniformly bounded above, and $\bar{\alpha}_g$ has a uniform positive lower bound, all independently of $n$, then $\|E_{\mathrm{TG}}\|_A\leq q<1$ for some constant $q$ independent of $n$.
\end{theorem}

Under the stated assumptions, the approximation constants $\bar{\beta}_g$ are uniformly bounded above independently of $n$.
For example, $\omega$-Jacobi smoothing with fixed $\omega\in(0,1)$ satisfies $\bar{\alpha}_g\geq[1-(2\omega-1)^2]/2>0$, providing the required uniform positive lower bound.

Note that \cref{lem:convergence} applies to large sparse matrices arising from PDE discretizations.
For dense matrices, the theorem does not yield a convergence bound independent of the matrix size $n$.

\section{Numerical experiments}
\label{sec:numerics}
The proposed \CFUA{} method is integrated into JXPAMG~\cite{2023-JXPAMG}.
In this section, we demonstrate the effectiveness of \CFUA{} through numerical experiments on three classes of 3D problems: Poisson-type model equations, radiation diffusion equations in radiation hydrodynamics, and incompressible Navier--Stokes equations in complex flows.
In these experiments, we compare \CFCA{} with four widely used aggregation algorithms: \Pone{}~\cite{1996-SAMG,1995-N-UA-AMG}, \PW{}~\cite{2006-UA-PW,2012-UA-CD}, \MIS{}~\cite{2000-par-sa,2013-GPU-ParUA}, and \CWM{}~\cite{2013-CWM-UA,2021-parallel-CWM}.

\subsection{Experimental setup}
All experiments are performed in an MPI environment on a Sugon supercomputer, using a single compute node with 64 CPU cores.
UA-AMG is used as a preconditioner for FGMRES in all experiments.
The convergence tolerance is set to $10^{-8}$ unless otherwise stated.
The parameters are listed in \cref{tab:parameters-table}.

\begin{table}[!htbp]
    \centering
    \caption{Algorithm parameters.}
    \label{tab:parameters-table}
    \small\renewcommand{\arraystretch}{1.15}
    \setlength{\tabcolsep}{6pt}
    \begin{tabular}{cccccc}
        \toprule
        \makecell[c]{Restart\\parameter} & \makecell[c]{Strength\\threshold} & \makecell[c]{Coarsest-grid\\size limit} & \makecell[c]{Smoother} & \makecell[c]{Cycle \\ type} & \makecell[c]{C/F splitting\\ strategy} \\
        \midrule
        30 & 0.25 & 50 & Hybrid SGS & W-cycle & PMIS \\
        \bottomrule
    \end{tabular}
\end{table}

The W-cycle provides better convergence efficiency than the V-cycle in these tests and is therefore used in this paper.
In addition, PMIS offers an advantage in CPU time over coarsening algorithms such as HMIS and Falgout in these tests and is therefore used for C/F splitting in \CFCA{} in this paper.

\subsection{Model problems: 3D Poisson-type equations}
We consider the 3D Poisson-type problem
\begin{equation}
    \label{tab:m-matrix}
    \begin{cases}
        -\nabla\cdot(\kappa\nabla u) = f, & x\in\Omega, \\
        u = 0, & x\in\partial\Omega,
    \end{cases}
\end{equation}
where $\Omega=(0,1)^3$.
The problem is discretized using a seven-point finite difference stencil on a structured $128\times128\times128$ grid.
Four 3D diffusion test cases with different diffusion coefficients $\kappa$ are considered:
\begin{enumerate}
    \item Constant coefficient (Diff3D-Const): $\kappa=1$.
    \item Anisotropic (Diff3D-Ani): $\kappa=\operatorname{diag}(1,10^3,10^3)$.
    \item Discontinuous coefficient (Diff3D-Dis): $\kappa=10^3$ for $x\in[0.25,0.75]^3$ and $\kappa=1$ elsewhere.
    \item Random coefficient (Diff3D-Rand): $\kappa=(10^3)^\delta$, where $\delta\in[0,1]$ is randomly distributed over $\Omega$.
\end{enumerate}
Matrix information for the four cases is listed in \cref{tab:poisson-info}.

\begin{table}[!htbp]
    \centering
    \caption{Test case information for the 3D Poisson-type problems.}
    \label{tab:poisson-info}
    \small\renewcommand{\arraystretch}{1.15}
    \setlength{\tabcolsep}{6pt}
    \begin{tabular}{cccc}
        \toprule
        \makecell[c]{Case} & \makecell[c]{Coefficient type} & \makecell[c]{Matrix size\\$N$} & \makecell[c]{Nonzeros\\$\mathit{nnz}$} \\
        \midrule
        Diff3D-Const & Constant & \multirow{4}{*}{2,097,152} & \multirow{4}{*}{14,581,760} \\
        Diff3D-Ani & Anisotropic & & \\
        Diff3D-Dis & Discontinuous & & \\
        Diff3D-Rand & Random & & \\
        \bottomrule
    \end{tabular}
\end{table}

\begin{figure}[!htbp]
    \centering
    \includegraphics[width=0.70\textwidth]{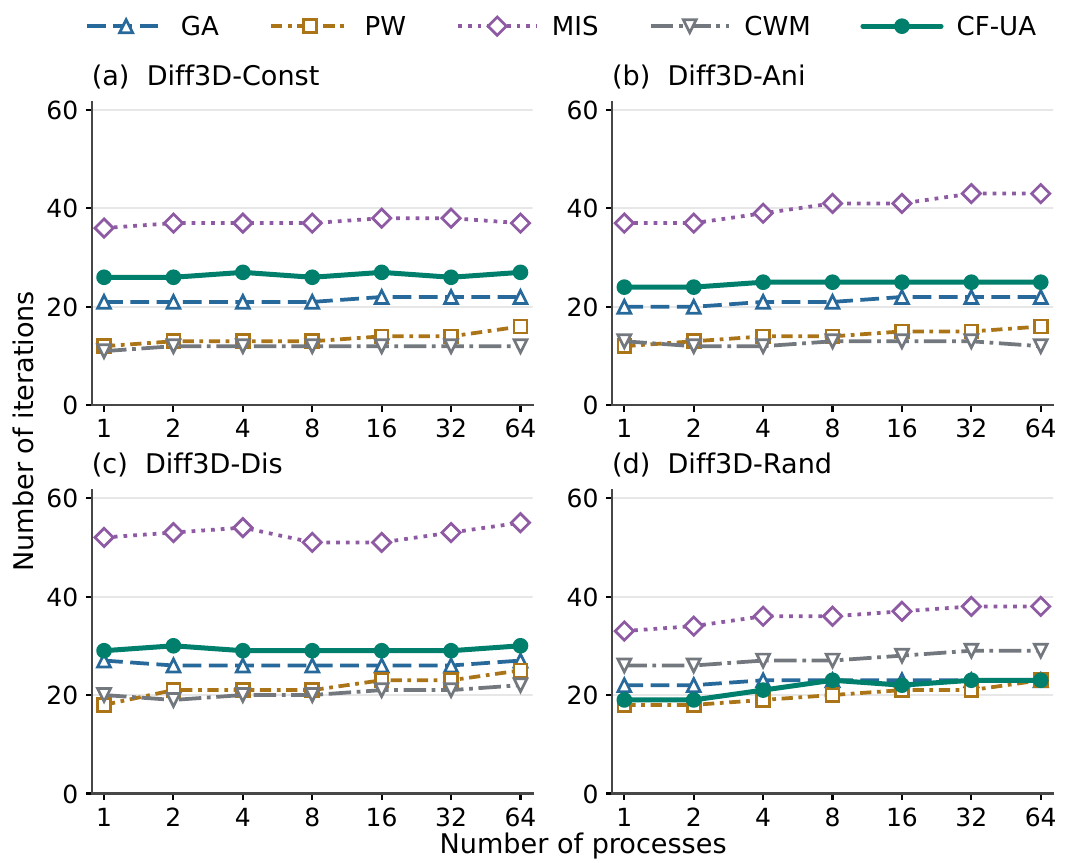}
    \caption{Number of iterations with different aggregation algorithms for the 3D Poisson-type problems.}
    \label{fig:poisson-iter}
\end{figure}

The variation in the number of iterations ($n_{\mathrm{it}}$) with the number of processes ($p$) is shown in \cref{fig:poisson-iter}.
For Diff3D-Const, Diff3D-Ani, and Diff3D-Dis, $n_{\mathrm{it}}$ for \CFCA{} remains nearly constant as $p$ increases, with values higher than those of \Pone{}, \PW{}, and \CWM{}.
For Diff3D-Rand, $n_{\mathrm{it}}$ for \CFCA{} shows only a modest increase with $p$ and remains lower than for \MIS{} and \CWM{}, with values comparable to those of \Pone{} and \PW{}.

\begin{figure}[!htbp]
    \centering
    \includegraphics[width=0.70\textwidth]{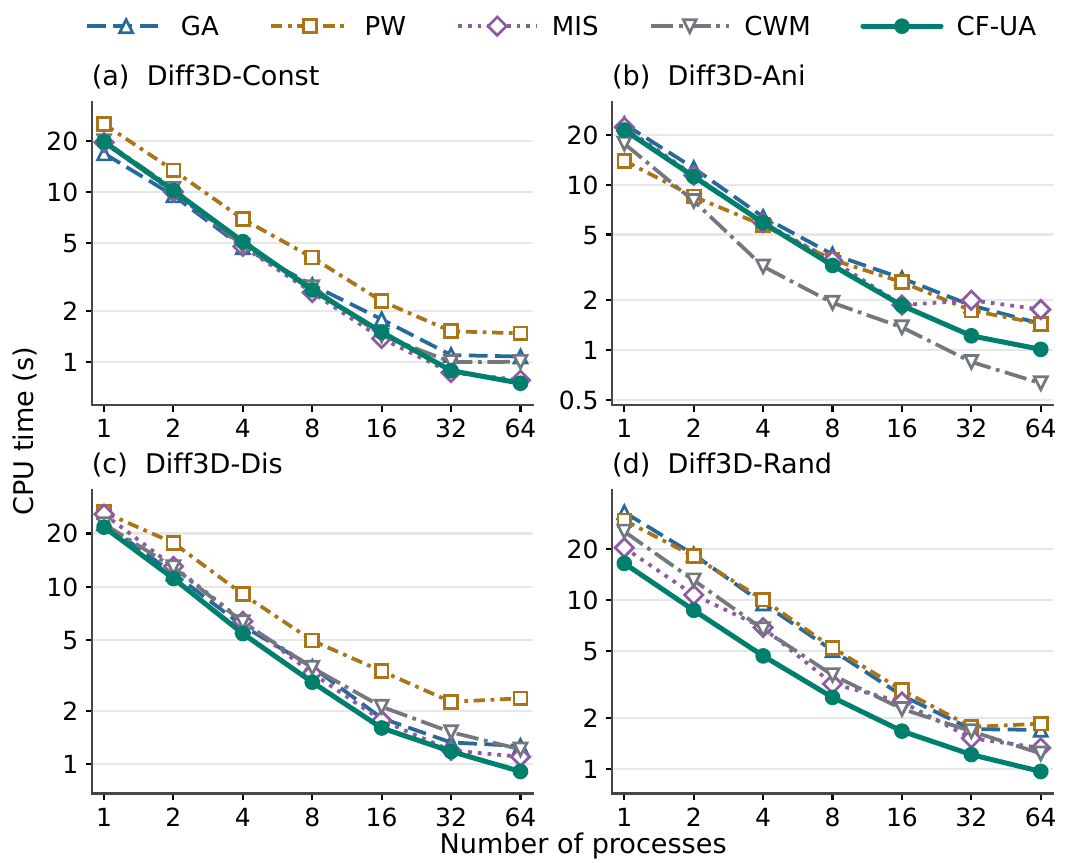}
    \caption{CPU times with different aggregation algorithms for the 3D Poisson-type problems.}
    \label{fig:poisson-time}
\end{figure}

CPU times for the four Poisson-type cases are compared in \cref{fig:poisson-time}.
For Diff3D-Const, the speedups of \CFCA{} relative to \Pone{}, \PW{}, \MIS{}, and \CWM{} over the tested values of $p$ are $0.85$--$1.44\times$, $1.27$--$1.97\times$, $0.92$--$1.04\times$, and $0.95$--$1.34\times$, respectively.
For Diff3D-Ani, the corresponding speedups are $1.08$--$1.53\times$, $0.65$--$1.43\times$, $1.00$--$1.74\times$, and $0.54$--$0.83\times$.
The advantage of \CFCA{} in CPU time is limited for the constant-coefficient and anisotropic cases. In particular, \CWM{} requires less CPU time than \CFCA{} for Diff3D-Ani at every tested value of $p$.
In contrast, for Diff3D-Dis, the corresponding speedups over the tested values of $p$ are $1.05$--$1.39\times$, $1.21$--$2.58\times$, $1.01$--$1.21\times$, and $1.04$--$1.33\times$.
For Diff3D-Rand, the corresponding speedups are $1.41$--$2.14\times$, $1.46$--$2.15\times$, $1.20$--$1.49\times$, and $1.28$--$1.54\times$.
For Diff3D-Dis and Diff3D-Rand, \CFCA{} achieves the lowest CPU times at every tested value of $p$, although it does not always yield the lowest $n_{\mathrm{it}}$ (see \cref{fig:poisson-iter}).
The reductions in CPU time are greater and more consistent for these two cases than for the constant-coefficient and anisotropic cases, demonstrating that \CFCA{} combines stable convergence with computational efficiency.

\subsection{Radiation hydrodynamics: 3D radiation diffusion equations}
We consider 3D radiation diffusion equations, which form an essential component of radiation hydrodynamics models for applications such as laser fusion.
The coupled energy equations for the radiation, electron, and ion temperatures are given by
\begin{equation}
    \left\{
    \begin{aligned}
        c_{vr}\frac{\partial T_r}{\partial t}
        - \frac{1}{\rho}\nabla\cdot(K_r\nabla T_r)
        &= \omega_{er}(T_e-T_r), \\
        c_{ve}\frac{\partial T_e}{\partial t}
        - \frac{1}{\rho}\nabla\cdot(K_e\nabla T_e)
        &= \omega_{ei}(T_i-T_e)+\omega_{er}(T_r-T_e), \\
        c_{vi}\frac{\partial T_i}{\partial t}
        - \frac{1}{\rho}\nabla\cdot(K_i\nabla T_i)
        &= \omega_{ei}(T_e-T_i).
    \end{aligned}
    \right.
\end{equation}
The system describes radiation diffusion and electron and ion heat conduction in a high-energy-density plasma.
These processes are coupled through electron--radiation and electron--ion energy exchange.
Here, $\rho$ is the material density, and $T_r$, $T_e$, and $T_i$ are the radiation, electron, and ion temperatures, respectively.
The coefficients $c_{vr}$, $c_{ve}$, and $c_{vi}$ are the corresponding effective specific heat capacities, and $K_r$, $K_e$, and $K_i$ are the corresponding heat-conduction (diffusion) coefficients.
The coefficients $\omega_{er}$ and $\omega_{ei}$ are the electron--radiation and electron--ion energy exchange coefficients, respectively.

This system is commonly referred to as the three-temperature (3T) model, with separate radiation, electron, and ion temperatures.
When radiation and matter are in energy equilibrium, the three temperatures coincide and the system reduces to the one-temperature (1T) model.
We consider one 1T case (RHD-1T) and three 3T cases (RHD-3T-1-1, RHD-3T-3-3, and RHD-3T-43-1).
The three 3T matrices correspond to different nonlinear iterations at different time steps.
The matrix characteristics are summarized in \cref{tab:rhd-info}.

\begin{table}[!htbp]
    \centering
    \caption{Test case information for the 3D radiation diffusion problems.}
    \label{tab:rhd-info}
    \small\renewcommand{\arraystretch}{1.15}
    \setlength{\tabcolsep}{6pt}
    \begin{tabular}{ccccc}
        \toprule
        \makecell[c]{Case} & \makecell[c]{Model} & \makecell[c]{Grid size} & \makecell[c]{Matrix size\\$N$} & \makecell[c]{Nonzeros\\$\mathit{nnz}$} \\
        \midrule
        RHD-1T & 1T & 2,097,152 & 2,097,152 & 14,581,760 \\
        \midrule
        RHD-3T-1-1 & \multirow{3}{*}{3T} & \multirow{3}{*}{2,097,152} & \multirow{3}{*}{6,291,456} & \multirow{3}{*}{52,133,888} \\
        RHD-3T-3-3 & & & & \\
        RHD-3T-43-1 & & & & \\
        \bottomrule
    \end{tabular}
\end{table}

\begin{figure}[!htbp]
    \centering
    \includegraphics[width=0.75\textwidth]{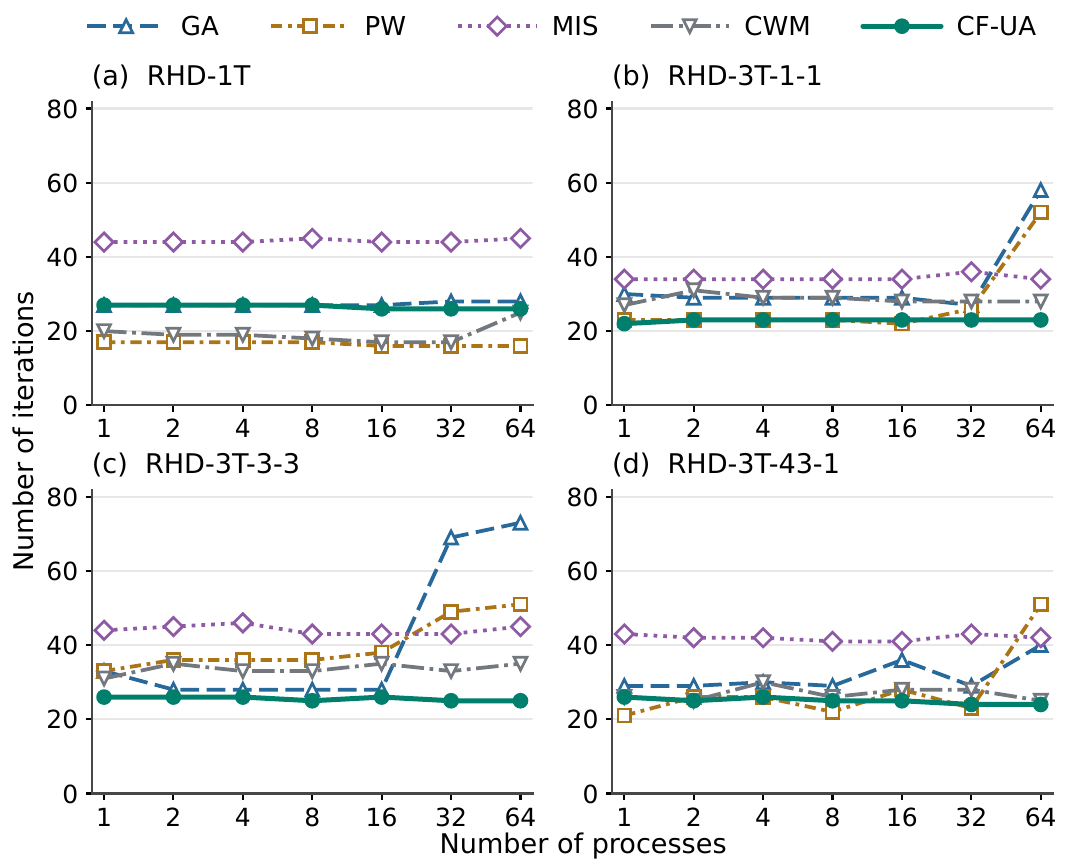}
    \caption{Number of iterations with different aggregation algorithms for the 3D radiation diffusion problems.}
    \label{fig:rhd-iter}
\end{figure}

The variation in $n_{\mathrm{it}}$ with $p$ for the four radiation diffusion cases is shown in \cref{fig:rhd-iter}.
For RHD-1T, $n_{\mathrm{it}}$ for \CFCA{} is comparable to that of \Pone{}, lower than that of \MIS{}, and higher than those of \PW{} and \CWM{}.
For the 3T cases, \CFCA{} maintains low and nearly constant values of $n_{\mathrm{it}}$ as $p$ increases.
In contrast, $n_{\mathrm{it}}$ for \Pone{} and \PW{} increases markedly at larger $p$ in some 3T cases, most notably RHD-3T-1-1 and RHD-3T-3-3.
The values for \MIS{} and \CWM{} also vary with $p$, with more pronounced variation for \CWM{} on RHD-3T-43-1.

\begin{figure}[!htbp]
    \centering
    \includegraphics[width=0.75\textwidth]{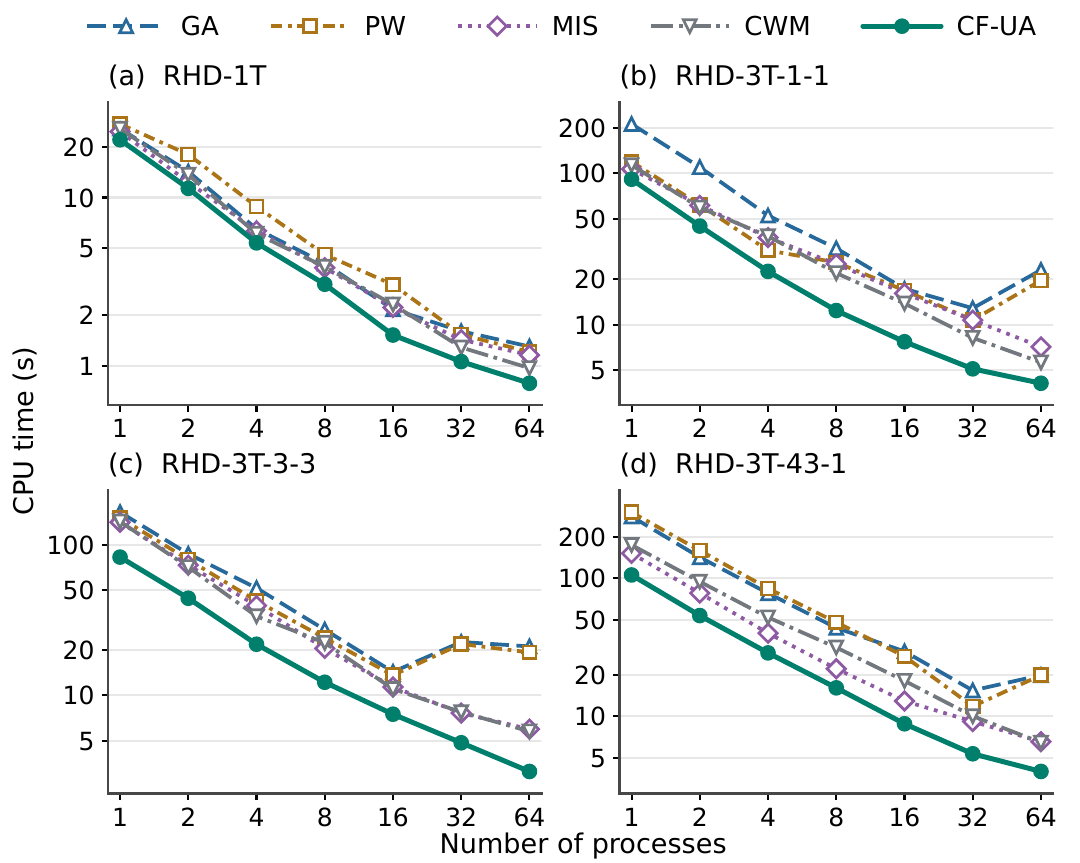}
    \caption{CPU times with different aggregation algorithms for the 3D radiation diffusion problems.}
    \label{fig:rhd-time}
\end{figure}

CPU times for the four radiation diffusion cases are compared in \cref{fig:rhd-time}.
\CFCA{} achieves the lowest CPU time for all four cases at every tested value of $p$.
For RHD-1T, its speedups relative to \Pone{}, \PW{}, \MIS{}, and \CWM{} over the tested values of $p$ are $1.16$--$1.65\times$, $1.24$--$1.99\times$, $1.09$--$1.47\times$, and $1.13$--$1.52\times$, respectively.
For RHD-3T-1-1, the corresponding speedups are $2.24$--$5.60\times$, $1.30$--$4.78\times$, $1.17$--$2.10\times$, and $1.24$--$1.79\times$.
For RHD-3T-3-3, the corresponding speedups are $1.91$--$6.75\times$, $1.80$--$6.16\times$, $1.52$--$1.91\times$, and $1.47$--$1.85\times$.
For RHD-3T-43-1, the corresponding speedups are $2.64$--$4.98\times$, $2.21$--$4.98\times$, $1.37$--$1.71\times$, and $1.62$--$2.05\times$.
The reduction in CPU time is particularly pronounced for the 3T cases, where \CFCA{} maintains low $n_{\mathrm{it}}$ as $p$ increases (see \cref{fig:rhd-iter}).
For RHD-3T-3-3, for example, \CFCA{} achieves speedups of up to $6.75\times$, $6.16\times$, $1.91\times$, and $1.85\times$ relative to \Pone{}, \PW{}, \MIS{}, and \CWM{}, respectively, with all four maxima attained at $p=64$.
These results demonstrate that \CFCA{} maintains stable convergence and computational efficiency for these complex problems.

\subsection{Complex flows: 3D incompressible Navier--Stokes equations}
We consider linear systems obtained by discretizing the 3D incompressible Navier--Stokes equations, the fundamental governing equations for many real-world complex flows:
\begin{equation}
    \label{IncompressibleNewtonianFluidNS}
    \left\{
    \begin{aligned}
        &\rho\frac{\partial\boldsymbol{u}}{\partial t}
        + \rho(\boldsymbol{u}\cdot\nabla)\boldsymbol{u}
        = -\nabla p + \mu\nabla^2\boldsymbol{u} + \boldsymbol{f}_b, \\
        &\nabla\cdot\boldsymbol{u} = 0,
    \end{aligned}
    \right.
\end{equation}
where $\boldsymbol{u}=(u,v,w)^T$ is the velocity vector, $p$ is the pressure, $\rho$ is the density, and $\mu$ is the dynamic viscosity.
Segregated and coupled methods are two common approaches to solving the incompressible Navier--Stokes equations.
For SIMPLE-type segregated methods, solving the pressure correction equation constitutes the main computational bottleneck~\cite{1999-SIMPLE,2002-book-CFD}.
For coupled methods, the main computational challenge is solving the coupled velocity--pressure linear system.

We consider two classes of test cases: lid-driven flows and cases from real-world applications.
All problems are discretized using a node-centered finite volume method~\cite{2016-book-FVM}.
We use nodal UA-AMG~\cite{2003-AMG-SYSTEM} to solve the coupled pressure--velocity systems, using a pressure-based measure matrix for coarsening~\cite{1996-N-UA-AMG}.

\subsubsection{Lid-driven flows}
For the pressure correction equation, we consider two domains: a unit cube (Pre3D-Cube) and a conical frustum (Pre3D-Cyl).
For the coupled pressure--velocity equations, we consider a box-shaped domain (Lid-3D).
All three cases are discretized on unstructured meshes, with the Reynolds number set to 1000.
The three cases are summarized in \cref{tab:test-lid}.

\begin{table}[!htbp]
    \centering
    \caption{Lid-driven flow test cases for the 3D incompressible Navier--Stokes equations.}
    \label{tab:test-lid}
    \small\renewcommand{\arraystretch}{1.15}
    \setlength{\tabcolsep}{6pt}
    \begin{tabular}{cccc}
        \toprule
        \makecell[c]{Case} & \makecell[c]{Type} & \makecell[c]{Grid size} & \makecell[c]{Matrix size} \\
        \midrule
        Pre3D-Cube & \multirow{2}{*}{Pressure correction} & 2,097,152 & 2,097,152 \\
        Pre3D-Cyl & & 6,403,951 & 6,403,951 \\
        \midrule
        Lid-3D & \makecell[c]{Coupled\\pressure--velocity} & 1,514,071 & 6,056,284 \\
        \bottomrule
    \end{tabular}
\end{table}

\begin{figure}[!htbp]
    \centering
    \includegraphics[width=1.0\textwidth]{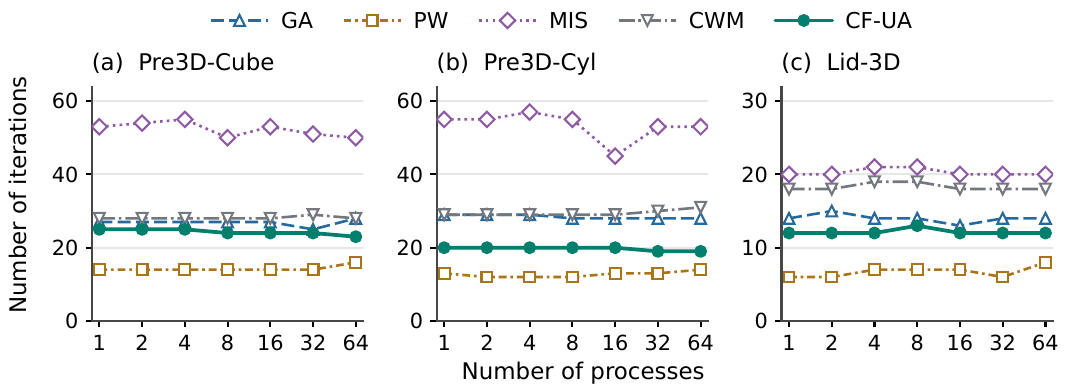}
    \caption{Number of iterations with different aggregation algorithms for the lid-driven flow problems.}
    \label{fig:lid-iter}
\end{figure}

The variation in $n_{\mathrm{it}}$ with $p$ for the three lid-driven flow cases is shown in \cref{fig:lid-iter}.
Across all three cases, \CFCA{} maintains nearly constant $n_{\mathrm{it}}$ as $p$ increases.
The values for \CFCA{} are consistently lower than those for \Pone{}, \MIS{}, and \CWM{}, whereas \PW{} yields the lowest $n_{\mathrm{it}}$.
This behavior demonstrates stable convergence for both pressure correction and coupled pressure--velocity systems.

\begin{figure}[!htbp]
    \centering
    \includegraphics[width=1.0\textwidth]{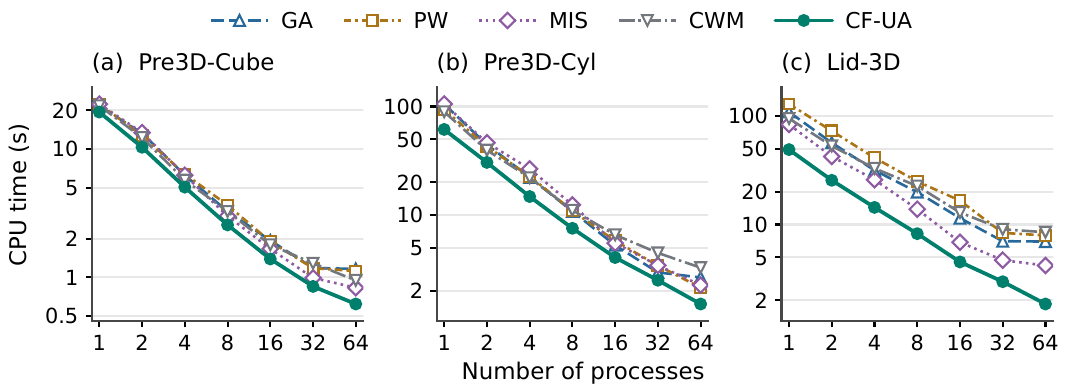}
    \caption{CPU times with different aggregation algorithms for the lid-driven flow problems.}
    \label{fig:lid-time}
\end{figure}

CPU times for the three lid-driven flow cases are compared in \cref{fig:lid-time}.
\CFCA{} achieves the lowest CPU time for all three cases at every tested value of $p$.
For Pre3D-Cube, its speedups relative to \Pone{}, \PW{}, \MIS{}, and \CWM{} over the tested values of $p$ are $1.12$--$1.88\times$, $1.16$--$1.81\times$, $1.15$--$1.34\times$, and $1.12$--$1.53\times$, respectively.
For Pre3D-Cyl, the corresponding speedups are $1.18$--$1.76\times$, $1.37$--$1.54\times$, $1.36$--$1.80\times$, and $1.30$--$2.16\times$.
For Lid-3D, the corresponding speedups are $2.20$--$3.76\times$, $2.63$--$4.29\times$, $1.52$--$2.26\times$, and $1.93$--$4.59\times$.
The reduction in CPU time is greater for the coupled Lid-3D problem: the speedups over \Pone{}, \PW{}, \MIS{}, and \CWM{} reach $3.76\times$, $4.29\times$, $2.26\times$, and $4.59\times$, respectively, at $p=64$.
\PW{} requires more CPU time than \CFCA{}, despite its lower $n_{\mathrm{it}}$ (see \cref{fig:lid-iter}).
Taken together, the results for $n_{\mathrm{it}}$ and CPU time demonstrate that \CFCA{} combines stable convergence with computational efficiency for these problems.

\subsubsection{Real-world cases}
We also consider two representative cases from real-world applications, each including 10 matrices from different time steps, all obtained by discretizing the coupled pressure--velocity equations.
The characteristics of these cases are summarized in \cref{tab:test-cases}.
The convergence tolerance for the relative pressure residual norm is set to $10^{-2}$ in both cases.

\begin{table}[!htbp]
    \centering
    \caption{Real-world cases for the 3D incompressible Navier--Stokes equations.}
    \label{tab:test-cases}
    \small\renewcommand{\arraystretch}{1.15}
    \setlength{\tabcolsep}{6pt}
    \begin{tabular}{ccccc}
        \toprule
        \makecell[c]{Case} & \makecell[c]{Number of\\matrices} & \makecell[c]{Type} & \makecell[c]{Grid size} & \makecell[c]{Matrix size} \\
        \midrule
        NS-C1 & \multirow{2}{*}{10} & \multirow{2}{*}{\makecell[c]{Coupled\\pressure--velocity}} & 60,075 & 240,300 \\
        NS-C2 & & & 591,000 & 2,364,000 \\
        \bottomrule
    \end{tabular}
\end{table}

\begin{figure}[!htbp]
    \centering
    \includegraphics[width=0.75\textwidth]{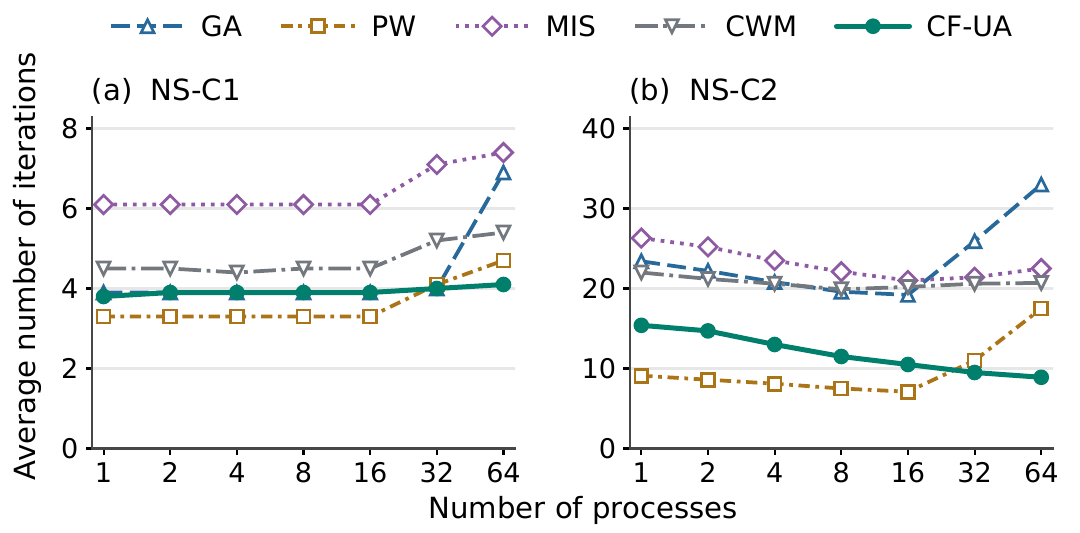}
    \caption{Average number of iterations with different aggregation algorithms for NS-C1 and NS-C2.}
    \label{fig:nsc12-iter}
\end{figure}

The variation in the average number of iterations ($\bar{n}_{\mathrm{it}}$) with $p$ is shown in \cref{fig:nsc12-iter}.
The average is taken over the 10 matrices in each case.
For NS-C1, \CFCA{} maintains an almost constant $\bar{n}_{\mathrm{it}}$ as $p$ increases, with values lower than those of \MIS{} and \CWM{}.
The values for \Pone{} and \PW{} are comparable to those for \CFCA{} at smaller $p$ and increase more noticeably at larger $p$.
For NS-C2, $\bar{n}_{\mathrm{it}}$ for \CFCA{} decreases steadily with $p$, whereas \Pone{} and \PW{} exhibit pronounced increases at larger $p$.
These results demonstrate the convergence robustness of \CFCA{} as $p$ increases for both real-world cases.

\begin{figure}[!htbp]
    \centering
    \includegraphics[width=1.0\textwidth]{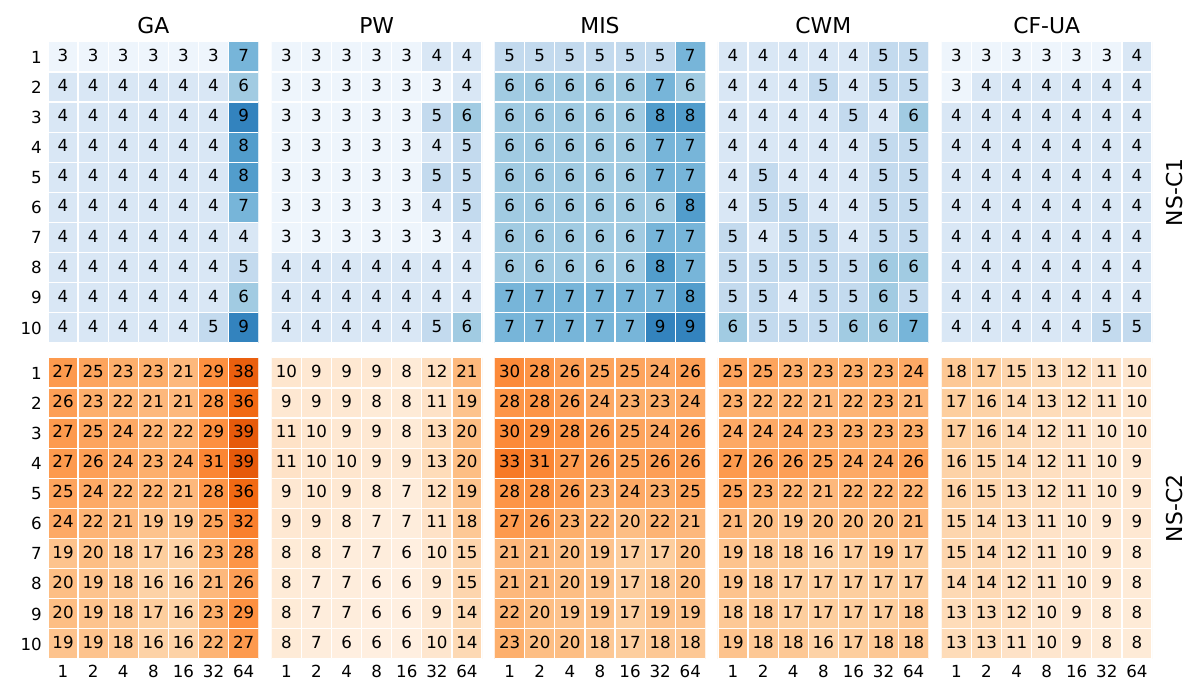}
    \caption{Number of iterations at individual time steps with different aggregation algorithms for NS-C1 and NS-C2.}
    \label{fig:nsc1-iter-all}
\end{figure}

The variation in $n_{\mathrm{it}}$ with $p$ at individual time steps is shown in the heatmaps in \cref{fig:nsc1-iter-all}.
The columns correspond to \Pone{}, \PW{}, \MIS{}, \CWM{}, and \CFCA{} from left to right, and the upper and lower rows correspond to NS-C1 and NS-C2, respectively.
In each panel, the horizontal and vertical axes indicate $p$ and the matrix index, respectively, and each cell gives $n_{\mathrm{it}}$ for that configuration.
Darker shading indicates higher $n_{\mathrm{it}}$, with separate color scales for the two cases and a common scale for all algorithms within each case.
For NS-C1, \CFCA{} maintains low $n_{\mathrm{it}}$ with limited variation across both $p$ and time steps, whereas \Pone{} and \PW{} show greater variation at the largest $p$.
For NS-C2, $n_{\mathrm{it}}$ for \CFCA{} decreases with $p$ and is nonincreasing as the time steps advance at each fixed $p$.
At the largest $p$, its values are lower than those of \Pone{} and \PW{}.
These results demonstrate the robustness of \CFCA{} with respect to both changes in matrix entries over time and increases in $p$.

\begin{figure}[!htbp]
    \centering
    \includegraphics[width=0.75\textwidth]{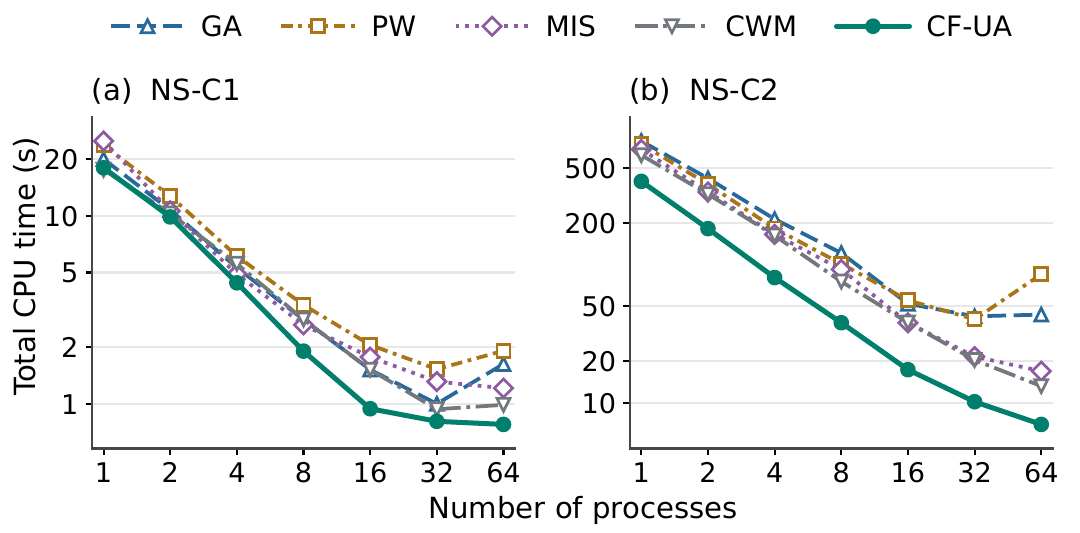}
    \caption{Total CPU times with different aggregation algorithms for NS-C1 and NS-C2.}
    \label{fig:nsc12-time}
\end{figure}

Total CPU times for solving all 10 matrices in each real-world case are compared in \cref{fig:nsc12-time}.
For NS-C1, the speedups of \CFCA{} relative to \Pone{}, \PW{}, \MIS{}, and \CWM{} over the tested values of $p$ are $1.09$--$2.10\times$, $1.29$--$2.47\times$, $1.08$--$1.88\times$, and $0.98$--$1.61\times$, respectively.
\CFCA{} achieves the lowest total CPU time at every tested value of $p$ except $p=1$, where \CWM{} requires less total CPU time.
For NS-C2, \CFCA{} achieves the lowest total CPU time at every tested value of $p$, with corresponding speedups of $1.96$--$6.23\times$, $1.86$--$12.20\times$, $1.70$--$2.43\times$, and $1.55$--$2.19\times$.

Total CPU time reductions are more pronounced for NS-C2, where \Pone{} and \PW{} exhibit a marked increase in $\bar{n}_{\mathrm{it}}$ at larger $p$, whereas \CFCA{} maintains low values (see \cref{fig:nsc12-iter}).
The speedups of \CFCA{} over \Pone{}, \PW{}, \MIS{}, and \CWM{} reach $6.23\times$, $12.20\times$, $2.43\times$, and $2.19\times$, respectively.
All four maximum speedups are attained with 64 processes.

\section{Conclusions}
\label{sec:conclusion}
A C/F splitting-based parallel UA-AMG method, \CFUA{}, has been proposed.
Uniform two-grid convergence has been established under the stated assumptions, with a convergence bound independent of the matrix size.
Numerical experiments on three classes of 3D problems demonstrate that \CFUA{} achieves lower CPU times in most comparisons with UA-AMG methods employing four widely used aggregation algorithms.
For the Poisson-type problems, the advantage in CPU time is more consistent for discontinuous and randomly varying coefficients than for the constant-coefficient and anisotropic cases.
The advantage of \CFUA{} in CPU time is greater for radiation diffusion and incompressible Navier--Stokes problems, particularly for the three-temperature cases and coupled pressure--velocity systems.
Across the tested problems, convergence remains robust as the number of processes increases.
This stability is also observed for the sequences of systems arising at different time steps in the real-world flow problems.

\section*{Acknowledgments}
The authors thank Prof.~Shi Shu, Prof.~Xiaoqiang Yue, and Dr.~Jianmeng He for their valuable advice on code design and their assistance with the implementation of the numerical methods used in this work.

\bibliographystyle{siamplain-citeorder}
\bibliography{references}

\end{document}